\documentclass[11pt,twoside,openright]{amsart}        

\usepackage{amsfonts}
\usepackage{amsmath}
\usepackage{amsthm}
\usepackage{amscd}
\usepackage{euscript}
\usepackage{amssymb}
\usepackage{psfrag}
\usepackage{graphicx}   
\usepackage{rotating}
\usepackage{color}
\usepackage{mathrsfs}
\usepackage{xypic}
\theoremstyle{plain}
\newtheorem{thm}{Theorem}[section]
\newtheorem{cor}[thm]{Corollary}
\newtheorem{lem}[thm]{Lemma}
\newtheorem{prop}[thm]{Proposition}
\theoremstyle{definition}
\newtheorem{defn}{Definition}[section]
\theoremstyle{remark}
\newtheorem{rem}{Remark}[section]
\theoremstyle{example}

\numberwithin{equation}{section}

\newcommand{\Z}{\mathbb Z}
\newcommand{\F}{\mathbb F}
\newcommand{\C}{\mathbb C}
\newcommand{\R}{\mathbb R}

\newcommand{\G}{\mathbb G}
\newcommand{\PP}{{\mathbb P}}

\newcommand{\OO}{\mathcal{O}}

\newcommand{\TT}{\mathcal{T}}

\newcommand{\Q}{\mathcal{Q}}

\DeclareMathOperator{\HH}{H} 
 \DeclareMathOperator{\Hom}{Hom}
 \DeclareMathOperator{\rk}{rk}

 \DeclareMathOperator{\Id}{Id}
\DeclareMathOperator{\T}{T} \DeclareMathOperator{\Spec}{Spec}
\DeclareMathOperator{\Proj}{Proj}

 \DeclareMathOperator{\Lie}{Lie}
\DeclareMathOperator{\Pic}{Pic} \DeclareMathOperator{\GL}{GL}
\DeclareMathOperator{\End}{End} 
 \DeclareMathOperator{\SL}{SL}

 \DeclareMathOperator{\gr}{gr}
\DeclareMathOperator{\SO}{SO}

\DeclareMathOperator{\lieg}{\mathfrak{g}}
\DeclareMathOperator{\liep}{\mathfrak{p}}

\DeclareMathOperator{\lieh}{\mathfrak{h}}
\DeclareMathOperator{\lier}{\mathfrak{r}}
\DeclareMathOperator{\lien}{\mathfrak{n}}

\title{On simplicity and stability of tangent bundles of rational homogeneous varieties}

\subjclass[2000]{14F05,14M17,14D20,16G20}
\keywords{Homogeneous Vector Bundle, Simplicity, Stability, Quiver Representation}

\author{Ada Boralevi}
\address{%
Department of Mathematics\\
Mailstop 3368\\
Texas A$\&$M University\\
College Station, TX 77843-3368\\
	USA}
\email{boralevi@math.tamu.edu}

\date{2009/01/10}

\begin{document}

\maketitle

\begin{abstract}
Given a rational homogeneous variety $G/P$ where $G$ is complex simple and of type $ADE$, we prove that all tangent bundles $\T_{G/P}$ are simple, meaning that their only endomorphisms are scalar multiples of the identity. This result combined with Hitchin-Kobayashi correspondence implies stability of these tangent bundles with respect to the anticanonical polarization. Our main tool is the equivalence of categories between homogeneous vector bundles on $G/P$ and finite dimensional representations of a given quiver with relations.
\end{abstract}

\section{Introduction}
In \cite{Rama} Ramanan proved that irreducible homogeneous bundles on rational homogeneous varieties are stable, and hence in particular simple. If the underlying variety is Hermitian symmetric then this result applies to tangent bundles. For the general case, the Hitchin-Kobayashi correspondence gives a weaker result for the tangent bundles, polystability. In this paper we show that in fact tangent bundles of any $G/P$ are simple, where $G$ is complex, simple and of type $ADE$. Simplicity and polystability combined give stability.\\ 
Our main tool is the equivalence of categories between homogeneous bundles on $G/P$ and finite dimensional representations of a given quiver with relations. Once the machinery of this equivalence of categories is set up, the simplicity of tangent bundles turns out to be an immediate and surprisingly easy consequence of it. Indeed one only needs to look at endomorphisms of the bundle as endomorphisms of the associated quiver representations. 

Homogeneous vector bundles have been classically studied using another equivalence of categories, namely that between homogeneous bundles on $G/P$ and finite dimensional representations of the parabolic subgroup $P$. 

In \cite{BK} Bondal and Kapranov had the idea of associating to any rational homogeneous variety a quiver with relations. By putting the 
appropriate relations one gets the aforementioned equivalence of categories between $G$-homogeneous vector bundles on $G/P$ and 
finite dimensional representations of the quiver. The relations were later refined by Hille in \cite{Hi1}.

In \cite{ACGP} \'{A}lvarez-C\'{o}nsul and Garc\'{i}a-Prada gave an equivalent construction, while in \cite{OR} Ottaviani and Rubei used the quiver for computing cohomology, obtaining a generalization of the well-known Borel-Weil-Bott theorem holding on Hermitian symmetric varieties of $ADE$ type.

We describe both equivalence of categories and give details on the quiver, its relations and its representations in Section 2 and 3.\\

Section 4 and 5 contain results on simplicity and stability. We use the quiver to prove that homogeneous vector bundles whose associated quiver representation has a particular 
configuration---we call such bundles multiplicity free---are weakly simple, which means that their only $G$-invariant endomorphisms are scalar multiple of the identity. 
Our result holds on any $G/P$, where $G$ is a simple group of type $ADE$:\\

\textbf{Proposition A.} \emph{Let $E$ be a multiplicity free homogeneous vector bundle of rank $r$ on $G/P$. 
Let $k$ be the number of connected components of the quiver $\Q|_E$. Then $\HH^0(\End E)^\C=\C^k$. In particular if $\Q|_E$ is connected, then $E$ is weakly simple.}\\

It turns out that all tangent bundles $\T_{G/P}$ are multiplicity free and connected, and that moreover the isotypical component 
$\HH^0(\End \T_{G/P})^\C$ coincides with the whole space $\HH^0(\End \T_{G/P})$, or in other words that these bundles are simple.\\

\textbf{Theorem B.} \emph{Let $\T_{G/P}$ the tangent bundle on a rational homogeneous variety $G/P$, where
$G$ is a complex simple Lie group of type $ADE$ and $P$ one of its parabolic subgroups. Then $\T_{G/P}$ is simple.}\\

If algebraic geometry, representation theory and quiver representations give us simplicity, for stability differential geometry also joins the team. A homogeneous variety $G/P$ is in fact also a K\"{a}hler manifold, and as such it admits a Hermite-Einstein structure. In virtue of the Hitchin-Kobayashi correspondence this is equivalent to the polistability of its tangent bundle. This together with Theorem B gives:\\

\textbf{Theorem C.} \emph{Let $\T_{G/P}$ the tangent bundle on a rational homogeneous variety $G/P$, where
$G$ is a complex simple Lie group of type $ADE$, and $P$ one of its parabolic subgroups.
Then $\T_{G/P}$ is stable with respect to the anticanonical polarization $-K_{G/P}$ induced by the
Hermite-Einstein structure.}\\

In the case where $G/P$ is a flag manifold point-hyperplane in $\PP^n$, we obtained a complete understanding of stability of the tangent bundle with respect to different polarizations:\\

\textbf{Proposition D.} \emph{Let $\F=\F(0,n-1,n)$ be the flag manifold of type $\SL_{n+1}/P(\alpha_1,\alpha_n)$, and set:
$$m(n)=\frac{-n+n\sqrt{4n^2+4n-3}}{2(n^2+n-1)}.$$
Then the tangent bundle $\T_\F$ is stable with respect to the polarization $\OO_\F(a,b)$ if and only if it is semistable if and only if $m(n)a \leq b \leq m(n)^{-1}a$.}\\

We also show similar computations for $\SL_4/B$.\\

In the last Section 6 we deal with moduli spaces. We quote and generalize the results from \cite{OR}, where the authors showed that
King's notion of semistability \cite{king} for a representation $[E]$ of the quiver $\Q_{G/P}$ is in fact equivalent to the
Mumford-Takemoto semistability of the associated bundle $E$ on $G/P$, when the latter is a Hermitian symmetric variety.\\
We can thus construct moduli spaces of $G$-homogeneous semistable bundles with fixed $\gr E$ on any homogenous variety $G/P$ of $ADE$ type.

\subsection*{Acknowledgements} This paper is part of my PhD thesis. I am very grateful to my advisor Professor 
Giorgio Ottaviani for the patience with which he followed this work very closely and 
for always transmitting me lots of encouragement and mathematical enthusiasm. 
I would also like to thank Professor Jean-Pierre Demailly for inviting me to Grenoble and for his warm hospitality.

\section{Preliminaries}

\subsection{Notations and first fundamental equivalence of categories}
Let $G$ be a complex semisimple Lie
group. We make a choice $\Delta=\{\alpha_1,\ldots,\alpha_n\}$ of simple roots of
$\lieg=\Lie G$, and we call $\Phi^+$ (respectively $\Phi^-$) the set of
positive (negative) roots. We denote by $\lieh \subset
\lieg$ the Cartan subalgebra so that $\lieg$ decomposes as:
$$\lieg=\lieh \oplus \bigoplus_{\alpha \in \Phi^+}\lieg_{\alpha} \oplus \bigoplus_{\alpha \in \Phi^-} \lieg_{\alpha}.$$
A parabolic subgroup $P \leq G$ is a subgroup $P=P(\Sigma)$, where:
$$\Lie(P(\Sigma))=\lieh \oplus \bigoplus_{\alpha \in \Phi^+}\lieg_{\alpha} \oplus \bigoplus_{\alpha \in \Phi^-(\Sigma)} \lieg_{\alpha},$$
for a subset $\Sigma \subset \Delta$ that induces $\Phi^-(\Sigma)=\{\alpha \in \Phi^- | \alpha=\sum_{\alpha_i \notin \Sigma}p_i\alpha_i\}$.\\
If $\Sigma=\Delta$, then $P(\Delta)=B$ is the Borel subgroup.

A rational homogeneous variety is a quotient $G/P$.

A vector bundle $E$ on $G/P$ is called ($G$)-homogeneous if there is an action of $G$ on $E$ such that the following diagram commutes:
$$\xymatrix{G \times E \ar[r]\ar[d]&E\ar[d]\\
G \times G/P \ar[r]& G/P}$$
where the bottom row is just the natural action of $G$ on the cosets $G/P$.\\
Note that the tangent bundle $\T_{G/P}$ on any rational homogeneous variety $G/P$ is obviously a $G$-homogeneous bundle.

The category of $G$-homogeneous vector bundle on $G/P$ is equivalent to the category $P$-mod of representations of $P$, and 
also to the category $\liep$-mod, where $\liep=\Lie P$, see for example \cite{BK}.

More in detail, the group $G$ is a principal bundle over $G/P$
with fiber $P$. Any $G$-homogeneous vector bundle $E$ of rank $r$ is induced by this principal bundle via a representation
$\rho:P\rightarrow \GL(r)$. We denote $E=E_\rho$. Indeed, $E$ of rank $r$ over $G/P$ is homogeneous
if and only if there exists a representation $\rho: P \rightarrow
\GL(r)$ such that $E \simeq E_{\rho}$, and this entails the aforementioned equivalence 
of categories.

For any weight $\lambda$ we denote by $E_\lambda$ the homogeneous bundle corresponding to the 
irreducible representation of $P$ with maximal weight $\lambda$. 
Here $\lambda$ belongs to the fundamental Weyl chamber of the reductive part of $P$. 
Indeed, $P$ decomposes as $P=R\cdot N$ into a reductive part $R$ and a nilpotent part $N$. 
At the level of Lie algebras this decomposition entails a splitting $\liep=\lier \oplus \lien$, with the obvious notation $\lier=\Lie R$ and $\lien=\Lie N$. 
Moreover from a result by Ise \cite{Ise} we learn that a representation of $\liep$ is completely reducible if and only if
it is trivial on $\lien$, hence it is completely determined by its restriction to $\lier$.\\
The well-known Borel-Weil-Bott theorem \cite{Bott} computes the cohomology of such $E_\lambda$'s by using purely Lie algebra tools, 
namely the action of the Weyl group on the weight $\lambda$. In particular the theorem states that if $\lambda$ is dominant then 
$\HH^0(E_\lambda)=\Sigma_\lambda$ (the irreducible representation of highest weight $\lambda$) and all higher cohomology vanishes.

\section{The quiver $\Q_{G/P}$}

\subsection{Definition of the quiver $\Q_{G/P}$ and its representations}\label{sezione def of the quiver}

Other than looking at homogeneous bundles as $P$-modules, it is useful to 
try a different point of view and look at these same bundles as representation of a given quiver with relations. 
For basics on quiver theory we refer the reader to \cite{DeWe} or \cite{king}.

To any rational homogeneous variety $G/P$ we can associate a quiver with relations, that we denote by $\Q_{G/P}$. 
The idea is to exploit all the information given by the choice of the parabolic subgroup $P$, with its decomposition $P=R\cdot N$. 

Let $\Lambda$ be the fundamental Weyl chamber of $G$, and let $\Lambda^+$ be the Weyl chamber of the reductive part $R$. Then we can give the following:
\begin{defn}\label{defQX}
Let $G/P$ be any rational homogeneous variety. The quiver $\Q_{G/P}$ is constructed as follows. The set of vertices is:
$$\Q_0=\Lambda^+=\{\lambda\:| \lambda\:\hbox{ dominant weight for }\:R\}.$$
There is an arrow connecting the vertices $\lambda$ and $\mu$ if and only if the vector space $\Hom(\lien \otimes E_\lambda,
E_\mu)^G$ is non-zero.
\end{defn}

\begin{rem}\label{quante frecce}
Definition \ref{defQX} is precisely the original one of Bondal and
Kapranov \cite{BK}, later used also by Alvarez-C\'{o}nsul and
Garc\'{\i}a-Prada \cite{ACGP}. Arrows correspond to weights of
the nilpotent algebra $\lien$, considered as an $\lier$-module with the adjoint action.\\
In fact one could obtain an equivalent theory by considering the
same vertices with a smaller number of arrows, i.e. by taking
only weights of the quotient $\lien/[\lien,\lien]$. This is for example the choice made by Hille (see \cite{Hi2} and \cite{Hi3}).
\end{rem}

Note that vertices $\lambda$ of $\Q_{G/P}$ correspond to irreducible homogeneous bundles $E_\lambda$ on $G/P$.

The relations on the quiver $\Q_{G/P}$ will be defined in Section \ref{sezione equivalenza categorie}.\\

Let now $E$ be an homogeneous vector bundle over $G/P$: we want to
associate to it a representation $[E]$ of the quiver $\Q_{G/P}$. The bundle $E$ comes with a filtration:
\begin{equation}\label{filtrazione}
0 \subset E_1 \subset E_2 \subset \ldots \subset E_k=E,
\end{equation}
where each $E_i/E_{i-1}$ is completely reducible. Of course the
filtration does not split in general. We define $\gr E=\oplus_i E_i/E_{i-1}$ for any filtration
(\ref{filtrazione}). The graded bundle $\gr E$ does not depend on
the filtration: in fact it is given by looking at our $\liep$-module $E$ as a module over
$\lier$, so that it decomposes as a direct sum of irreducibles:
\begin{equation}\label{grE}
\gr E = \bigoplus_\lambda E_\lambda \otimes V_\lambda,
\end{equation}
with multiplicity spaces $V_\lambda \simeq \C^k$, where $k \in \Z_{\geq 0}$
is the number of times $E_\lambda$ occurs.

The representation $[E]$ associates to the vertex $\lambda$ of the quiver 
$\Q_{G/P}$ precisely the multiplicity space $V_\lambda$ in the decomposition (\ref{grE}).

Going on with the definition of the representation $[E]$, given any $\lambda, \mu \in \Q_0$ such that there is an arrow $\lambda \rightarrow \mu \in \Q_1$, 
we now need to define a linear map $V_\lambda \rightarrow V_\mu$. This information is given by the nilpotent part $\lien$. 
More precisely, it is given by the natural action of $\lien$ over $\gr E$,
both viewed as $\lier$-modules:
$$\theta: \lien \otimes \gr E \rightarrow \gr E.$$
The morphism $\theta$ encodes all the information we need, including that on the relations of the quiver.

Obviously if we have a vector bundle $E$ we then have the graded $\gr E$ and the morphism $\theta$. 
Viceversa, if we have an $\lier$-module and a morphism that behaves ``just like $\theta$'', we can reconstruct a $\liep$-module and hence a vector bundle. 
More in detail, let us state and prove the following generalization of \cite[Theorem 3.1]{OR}:
\begin{thm}\label{theta} Consider $\lien$ as an $\lier$-module with the adjoint action.
\begin{enumerate}
  \item Given a $\liep$-module $E$ on $X$, the action of $\lien$ over $E$ induces a morphism of $\lier$-modules:
  $$\theta: \lien \otimes \gr E \rightarrow \gr E.$$
The morphism
$$\theta \wedge \theta:\wedge^2\lien\otimes \gr E \rightarrow \gr E $$
defined by $\theta \wedge \theta ((n_1\wedge n_2)\otimes f):= n_1 \cdot (n_2 \cdot f) -n_2 \cdot (n_1 \cdot f)$
satisfies the equality $\theta \wedge \theta=\theta \varphi$ in $\Hom(\wedge^2\lien\otimes \gr E, \gr E)$,
where $\varphi$ is given by:
  \begin{align}
  \nonumber \varphi:\wedge^2 \lien \otimes \gr E &\rightarrow \lien \otimes \gr E\\
  \nonumber (n_1 \wedge n_2)\otimes e &\mapsto [n_1,n_2] \otimes e.
  \end{align}
  \item Conversely, given an $\lier$-module $F$ on $X$ and a morphism of $\lier$-modules
  $$\theta: \lien \otimes F \rightarrow F$$
  such that $\theta \wedge \theta=\theta \varphi$, we have that $\theta$ extends uniquely to an action of $\liep$ over $F$,
  giving a bundle $E$ such that $\gr E=F$.
\end{enumerate}
\end{thm}
\begin{proof}
\begin{enumerate}
  \item Obviously $\theta$ is $\lier$-equivariant, almost by definition.\\
Take $n_1$ and $n_2$ in $\lien$. We have that:
  $$\theta \wedge \theta ((n_1\wedge n_2)\otimes f)= n_1 \cdot (n_2 \cdot f) -n_2 \cdot (n_1 \cdot f)=[n_1,n_2]\cdot f,$$
  which means exactly that $\theta \wedge \theta= \theta \varphi$.
  \item For any $f \in F$ and any $r+n \in \liep= \lier \oplus
  \lien$ we set:
  \begin{equation}\label{azione}
  (r+n)\cdot f := r \cdot f + \theta (n \otimes f).
  \end{equation}
We need to show that given any $p_1,p_2 \in \liep$, the action
(\ref{azione}) respects the bracket, i.e. that for every $f\in F$:
$$[p_1,p_2]\cdot f= p_1 \cdot (p_2 \cdot f) -p_2 \cdot (p_1 \cdot f).$$
Now if both $p_1,p_2 \in \lier$, then there is nothing to prove.\\
If both $p_1,p_2 \in \lien$, then from the equality $\theta
\wedge \theta= \theta \varphi$, we get:
\begin{align}
\nonumber [p_1,p_2]\cdot f&=\theta ([p_1,p_2] \otimes f)= \theta \varphi ((p_1 \wedge p_2) \otimes f)=\\
\nonumber &=\theta \wedge \theta((p_1 \wedge p_2) \otimes f)=p_1 \cdot (p_2 \cdot f) -p_2 \cdot (p_1 \cdot f). 
\end{align}
Finally, in case $p_1 \in \lier$ and $p_2 \in \lien$, then
$[p_1,p_2] \in \lien$ and we have:
\begin{align}
\nonumber p_1 \cdot (p_2 \cdot f)&=\theta (p_1 \cdot(p_2 \otimes f))=\theta ( ([p_1,p_2] \otimes f + p_2 \otimes (p_1 \cdot f))=\\
\nonumber &=[p_1,p_2] \cdot f + \theta (p_2 \otimes (p_1 \cdot f))=[p_1,p_2] \cdot f + p_2 \cdot (p_1 \cdot f).
\end{align}
\end{enumerate}
\end{proof}

Remark that (\ref{grE}) entails that we have a decomposition:
\begin{equation}\label{decomposizione theta}
\theta \in \Hom(\lien \otimes \gr E, \gr E)= \bigoplus_{\lambda, \mu
\in \Q_0} \Hom(V_\lambda,V_\mu) \otimes
 \Hom(\lien \otimes E_\lambda,E_\mu).
\end{equation}

\begin{lem}\label{key lemma}\cite[Proposition 2]{BK}
In the $ADE$ case $\dim \Hom(\lien \otimes E_\lambda, E_\mu)^G$ is
either 0 or 1 for every $\lambda, \mu \in \Lambda^+$.
\end{lem}

\begin{rem}
From now on $G$ will thus always denote a complex Lie group of $ADE$ type. Nevertheless, the construction of the quiver with its relation can be done
for any type of Lie group, like in \cite{ACGP}.
\end{rem}

We can now conclude the construction of the representation of the quiver $[E]$
associated to the bundle $E$.\\
For any $\lambda$ $\lier$-dominant weight fix a maximal vector $v_\lambda$ (it is unique up to
constants). For any $\alpha$ weight of $\lien$, fix an eigenvector $n_\alpha \in \lien$.

Now suppose that there is an arrow $\lambda \rightarrow \mu$ in the quiver. Then the vector space $\Hom(\lien \otimes
E_\lambda,E_\mu)^G$ is non-zero, and in particular is 1-dimensional. Notice that by definition, being given by the action of
$\lien$, the arrow will send the
weight $\lambda$ into a weight $\mu=\lambda+\alpha$, for some
$\alpha \in \Phi^+$ root of $\lien$
(for we have $\lieg_\alpha \cdot W_\lambda \subseteq W_{\lambda+\alpha}$).\\
Then fix the generator $f_{\lambda\mu}$ of $\Hom(\lien \otimes
E_\lambda,E_\mu)^G$ that takes $n_\alpha \otimes v_\lambda \mapsto v_\mu$.\\
Once all the generators are fixed, from (\ref{decomposizione theta}) we write the map $\theta$ uniquely as:
\begin{equation}\label{decomposizione theta bis}
	\theta= \sum_{\lambda, \mu}g_{\lambda\mu}f_{\lambda\mu},
\end{equation}
and thus we can associate to the arrow $\lambda
\xrightarrow{f_{\lambda\mu}} \mu$ exactly the element
$g_{\lambda \mu}$ in $\Hom(V_\lambda,V_ \mu)$. All in all:
\begin{defn}\label{def rappresentazione}
To any homogeneous vector bundle $E$ on $G/P$ we associate a representation $[E]$ of the quiver $\Q_{G/P}$ as follows.\\
To any vertex $\lambda \in \Q_0$ we associate the vector space $V_\lambda$ from the decomposition (\ref{grE}).\\
To any arrow $\lambda \rightarrow \mu$ we associate the element $g_{\lambda \mu} \in \Hom(V_\lambda,V_\mu)$ from the decomposition
(\ref{decomposizione theta}).
\end{defn}

Notice that a different choice of generators would have led to an equivalent construction.

\subsection{Second fundamental equivalence of categories}\label{sezione equivalenza categorie}

We introduce here the equivalence of categories between homogeneous 
vector bundles on $G/P$ and finite dimensional representations of the quiver $\Q_{G/P}$.\\
From Proposition \ref{theta} it is clear that by putting the appropriate relations on the quiver, namely the equality 
$\theta \wedge \theta=\theta \varphi$, we can get the following:
\begin{thm}\cite{BK, Hi1, ACGP}\label{BK}
Let $G/P$ a rational homogeneous variety of $ADE$ type. The category of finite dimensional representations of the Lie
algebra $\liep$ is equivalent to the category of finite dimensional
representations of the quiver $\Q_{G/P}$ with certain relations
$\mathcal{R}$, and it is equivalent to the category of $G$-homogeneous
bundles on $G/P$.
\end{thm}

We show here how one can derive the relations. For details, see \cite{ACGP}.

Let $\lambda, \mu, \nu \in \Q_0$. We start by defining the morphism $\phi_{\lambda\mu\nu}$:
$$\Hom(\lien \otimes E_\lambda,E_\mu)\otimes \Hom(\lien \otimes E_\mu,E_\nu) \xrightarrow{\phi_{\lambda\mu\nu}}
\Hom(\wedge^2 \lien \otimes E_\lambda,E_\nu)$$
by setting $\phi_{\lambda\mu\nu}(a \otimes a'):(n \wedge n')\otimes x \mapsto z$, where $a:n \otimes x \mapsto y$ and $a':n' \otimes y\mapsto z$.

Obviously then the image of the $G$-invariant part:
$$\phi_{\lambda\mu\nu}(\Hom(\lien\otimes E_\lambda,E_\mu)^G\otimes \Hom(\lien \otimes E_\mu,E_\nu)^G)
\subseteq \Hom(\wedge^2 \lien \otimes E_\lambda,E_\nu)^G.$$
In particular recall once the choice of constants has been made, there are fixed generators $f_{\lambda\mu}$, where
$f_{\lambda\mu}: n_\alpha \otimes v_\lambda \mapsto v_\mu$, and $\alpha=\lambda -\mu$. Then if we set $\beta=\mu -\nu$, all in all:
$$\phi_{\lambda\mu\nu}(f_{\lambda\mu} \otimes f_{\mu\nu}): (n_\alpha \wedge n_\beta) \otimes v_\lambda \mapsto v_\nu.$$

Now consider the natural morphism
$\wedge^2 \lien \rightarrow \lien$ sending $n \wedge n' \mapsto [n,n']$.\\
It induces a morphism $\phi_{\lambda\nu}$:
$$\Hom(\lien \otimes E_\lambda,E_\nu) \xrightarrow{\phi_{\lambda\nu}} \Hom(\wedge^2 \lien \otimes E_\lambda,E_\nu).$$

Once again it is clear that the invariant part $\phi_{\lambda\nu}(\Hom(\lien \otimes E_\lambda,E_\nu)^G)$ is
contained in $\Hom(\wedge^2 \lien \otimes E_\lambda,E_\nu)^G$.

Theorem \ref{theta} together with the splitting (\ref{decomposizione theta bis}) entail that we have an equality in $\Hom(\wedge^2 \lien \otimes E_\lambda,E_\nu)^G$:
\begin{equation}\label{relazione da espandere}
\sum_{\lambda,\nu}\sum_\mu\Big(\phi_{\lambda\mu\nu}(f_{\lambda\mu} \otimes f_{\mu\nu})(g_{\lambda\mu}g_{\mu\nu})
+\phi_{\lambda\nu}([f_{\lambda\mu},f_{\mu\nu}])g_{\lambda\nu}\Big)=0.
\end{equation}

Finally, let $\{c_{\lambda\nu}^k\}$ be a basis of the vector space
$\Hom(\wedge^2 \lien \otimes E_\lambda,E_\nu)^G$, for
$k=1,\ldots,m_{\lambda \nu}$, $m_{\lambda \nu}=
\dim(\Hom(\wedge^2\lien \otimes E_\lambda,E_\nu)^G)$. Expand:
\begin{align}
\nonumber \phi_{\lambda\mu\nu}(f_{\lambda\mu} \otimes f_{\mu\nu})
=\sum_{k=1}^{m_{\lambda \nu}}  x_{\lambda\mu\nu}^k c_{\lambda \nu}^k\\
\nonumber \phi_{\lambda\nu}([f_{\lambda\mu},f_{\mu\nu}])
=\sum_{k=1}^{m_{\lambda \nu}}y_{\lambda\nu}^k c_{\lambda\nu}^k,\end{align}
Then for every couple of vertices $\lambda,\nu$, equality (\ref{relazione da espandere})
gives us a system of $m_{\lambda\nu}$ equations that the maps $g_{\gamma\delta}$ satisfy:
\begin{equation}\label{Rk}
\mathcal{R}_k^{\lambda\nu}: \sum_{\mu \in \Q_0} x_{\lambda\mu\nu}^k g_{\lambda \mu}g_{\mu\nu}+ y_{\lambda\nu}^k g_{\lambda\nu}=0
\end{equation}
\begin{defn}\label{def relazioni}
We define the relations $\mathcal{R}$ on the quiver $\Q_{G/P}$ as the
ideal generated by all the equations (that with a slight abuse of
notation we
 keep calling $\mathcal{R}_k$):
$$\mathcal{R}_k^{\lambda\nu}: \sum_{\mu \in \Q_0} x_{\lambda\mu\nu}^k
f_{\lambda \mu}f_{\mu\nu}+ y_{\lambda\nu}^k f_{\lambda\nu}=0,$$ for
$k=1,\ldots,m_{\lambda \nu}$ and for any couple of weights
$\lambda,\nu \in \Q_0$.
\end{defn}

\section{Simplicity}

\subsection{Simplicity of multiplicity free bundles} 
In this section we work on rational homogeneous varieties $G/P$, where $G$ is complex, simple of type $ADE$.

Let $E$ be a rank $r$ homogeneous vector bundle on $G/B$, and let $[E]$ be the associated representation of
the quiver $\Q_{G/P}$. Denote by $\Q|_E$ the subquiver of $\Q_{G/P}$ given by all vertices where $[E]$ is non-zero and all arrows connecting any two such vertices.
Clearly the support of $\Q|_E$ has at most $r$ vertices. Notice also that the representation $[E]$ of $Q_{G/P}$ induces a representation 
of the subquiver $\Q|_E$.\\
Call the vertices of $\Q|_E$ $\{\lambda_1,\ldots,\lambda_n\}$, with $n\leq r$. Then the usual decomposition for the graded of $E$ can be written:
\begin{equation}\label{deco}
\gr E=\bigoplus_{i=1}^n E_{\lambda_i} \otimes V_i.
\end{equation}
\begin{defn}
A homogeneous vector bundle $E$ of rank $r$ is \emph{multiplicity free}
if $\dim V_i=1$ for every $i=1,\ldots,n$ in (\ref{deco}).
\end{defn}

So let now $E$ be multiplicity free. Just by looking at the associated quiver representation $[E]$, we will show that the only 
$G$-invariant endomorphisms that such a bundle can have are scalar multiple of the identity, i.e. that the isotypical component 
$\HH^0(\End E)^\C=\C$. If this holds we call the bundle weakly simple.
\begin{prop}\label{prop sui 1^r-type}
Let $E$ be a multiplicity free homogeneous vector bundle of rank $r$ on $G/P$. 
Let $k$ be the number of connected components of the quiver $\Q|_E$. Then $\HH^0(\End E)^\C=\C^k$. In particular if $\Q|_E$ is connected, then $E$ is weakly simple.
\end{prop}
\begin{proof}
Suppose first that $k=1$, so that the subquiver $\Q|_E$ is connected.\\
Any element $\varphi \in \HH^0(\End E)^\C$ is a $G$-invariant endomorphism $\varphi:E \rightarrow E$. In particular we can
look at $\varphi$ as a morphism $[E] \rightarrow [E]$ between representations
of the same quiver. This means that we can look at $\varphi$ as a family of morphisms $\{\varphi_i:V_i \rightarrow V_i\:|\:i+1,\ldots,n\}$.
The hypothesis that $E$ is multiplicity free entails that in particular each $\varphi_i=k_i \Id$ and hence $\varphi=(k_1,\ldots,k_n)$
is in fact just an element of $\C^r$.\\
By definition of morphism of quiver representations, every time that there is an arrow $\lambda_i \rightarrow \lambda_j$ in the quiver
$\Q|_E$, there is a commutative diagram:
$$\xymatrix{V_i \ar[d]_{k_i}\ar[r] &V_j\ar[d]^{k_j}\\
V_i \ar[r] &V_j}$$ Remark that the two horizontal arrows are the
same. This means that if we fix the first constant $k_1$
of $\varphi$ then all $\varphi$ is completely determined thanks to connectedness. 
So this proves that $\HH^0(\End E)^\C \leq \C$.\\
On the other hand, notice we have homotheties, hence $\C \subseteq \HH^0(\End E)^\C$, and the thesis follows for the case $k=1$.\\
The same argument applies for each connected component, and this completes the proof.
\end{proof}

\subsection{Simplicity of tangent bundles}
The results of the previous section can be applied to a ``special'' multiplicity free homogeneous bundle: 
the tangent bundle $\T_{G/P}$, with $G$ simple of $ADE$ type.

Let $\T_{G/P}$ be the tangent bundle on $G/P$. Recall that any parabolic subgroup $P=P(\Sigma)$ is given by a
subset $\Sigma \subseteq \Delta$ of simple roots. Define also the subset of negative roots $\Phi^-_P=\Phi^-\setminus\Phi^-(\Sigma)$. 
Notice that when $P=B$ is the Borel, $\Sigma=\Delta$ and $\Phi^-_B=\Phi^-$.\\
The bundle $\T_{G/P}$ is a homogeneous bundle of rank $r=\dim G/P
=|\Phi^-_P|$, whose weights are exactly the elements of $\Phi^-_P$.

\begin{rem}\label{rem sul pullback}
It is convenient for us to take into account all the weights
of the tangent bundle, and not only the highest weights. Obviously
in the case of the Borel it doesn't make any
difference. If $P$ is any other parabolic subgroup of $G$, this
means that instead of the tangent bundle $\T_{G/P}$ we are
considering its pull-back $\pi^* \T_ {G/P}$ via the (flat!)
projection
$$\pi:G/B \rightarrow G/P.$$
The obvious vanishing $R^i\pi_*\OO=0$ for $i>0$, together with the
Projection formula (see \cite{Har}, II.5) guarantee that:
$$\HH^0(G/B, \End (\pi^*\T_{G/P}))=\HH^0(G/B,\pi^*(\End \T_{G/P}))=\HH^0(G/P,\End \T_{G/P}).$$
Hence we are allowed to work on $\pi^* \T_{G/P}$ instead of
$\T_{G/P}$. To simplify the notation we write $\TT_{G/P}:=\pi^* \T_{G/P}$.\\
Notice that the representation associated to any $\TT_{G/P}$, $P \supseteq
B$, is a subrepresentation of that associated to $\TT_{G/B}=\T_{G/B}$.
\end{rem}

Now we make the easy but fundamental remark that for every
homogeneous variety $G/P$, the rank $\TT_{G/P}$ (that is, the
dimension of $G/P$) coincides with the number of weights of the
associated representation and with the number of vertices in the
quiver representation $[\TT_{G/P}]$. Hence they must all have
multiplicity one. Moreover, notice that in the $ADE$ case whenever $\alpha \neq -\beta$ we have $[\lieg_\alpha,\lieg_\beta]=\lieg_{\alpha+\beta}$. All in all:
\begin{thm}\label{i tangenti sono 1r}
For every homogeneous variety $G/P$ the bundle $\TT_{G/P}$ is multiplicity free and connected.
\end{thm}

We will now show that in the case of a tangent bundle $\TT_{G/P}$  all endomorphisms are $G$-invariant endomorphisms, or in other words that 
$$\HH^0(\End(\TT_{G/P}))^\C=\HH^0(\End(\TT_{G/P})).$$
\begin{lem}\label{struttura H0 graduato}
Let $G/P$ any homogeneous variety of type $ADE$. For any $G$-module
$W \neq \C$, $\HH^0(\gr\End(\TT_{G/P}))^W=0$.
\end{lem}

\begin{proof}
This is a direct computation.\\
We use here the notation for positive and negative roots $\Phi^+$ and $\Phi^-$ used for example in \cite{FH}, and we call $\{L_i\}$ the
standard basis of $\lieh^*$. In the $A_{N-1}$ case ($G=\SL_N$):
$$\Phi^+=\{L_i-L_j\}_{1\leq i<j\leq N}\:\:\:\:\:\:\:\:
\hbox{and}\:\:\:\:\:\:\:\:\Phi^-=\{L_i-L_j\}_{1\leq j<i\leq N},$$
and the fundamental Weyl chamber associated is the set:
\begin{equation}\label{fund weyl chamber A_n}
\Lambda=\{\sum a_iL_i\:|\:a_1\geq a_2\geq \ldots \geq a_N\}.
\end{equation}
For the $D_N$ case ($G=\SO_{2N}$):
$$\Phi^+=(\{L_i-L_j\}\cup\{L_i+L_j\})_{i<j}\:\:\:\:\:
\hbox{and}\:\:\:\:\:\Phi^-=(\{L_i-L_j\}\cup\{-L_i-L_j\})_{i>j},$$
and the fundamental Weyl chamber associated:
\begin{equation}\label{fund weyl chamber D_n}
\Lambda=\{\sum a_iL_i\:|\:a_1\geq \ldots \geq a_{N-1}\geq |a_N|\}.
\end{equation}
By the Borel-Weil-Bott theorem \cite{Bott}, 
an irreducible bundle $E_\nu$ has non-zero $\HH^0(E_\nu)$ exactly if and only if $\nu$ lies in the fundamental Weyl chamber.

For the case of $A_{N-1}$: let $\nu$ be an irreducible summand of $\End \TT_{\SL_N/P}$. Hence $\nu$ is of the form
$\alpha + \beta$ with $\alpha \in \Phi^+$ and $\beta \in \Phi^-$.\\
Practically speaking, $\alpha$ is a vector in $\Z^N$ having $1$ at
the $i-$th place, $-1$ at the $j-$th place and $0$ everywhere else,
with $1\leq i < j \leq N$; similarly $\beta$ has a $-1$ at the
$h-$th place, a $1$ at the $k-$th place and $0$
everywhere else, with again $1\leq h < k \leq N$. What does a sum $\nu=\alpha+\beta$ look like? Fix $\alpha$. 
Then only six possibilities can occur for $\beta$, namely 
$$\left.
\begin{array}{ll}
\hbox{either:}&h<k \leq i<j,\\
\hbox{or:}&h \leq i \leq k\leq j,\\
\hbox{or:}&h \leq i <j \leq k,\\
\hbox{or:}&i \leq h <k \leq j,\\
\hbox{or:}&i \leq h \leq j \leq k,\\
\hbox{or else:}&i <j \leq h <k.
\end{array}
\right.
$$
One can check directly that 
the condition (\ref{fund weyl chamber A_n}) is satisfied exactly when $\beta=-\alpha$. In this case $\nu=0$ and $\HH^0(E_0)=\C$, which is what we wanted.

Let us now move to the case $D_N$. Here the situation is complicated by the fact that we deal with more roots, and thus with more
possible combinations for the sum $\nu=\alpha+\beta$. Nevertheless the condition (\ref{fund weyl chamber D_n}) for $\nu$ to belong
to the fundamental Weyl chamber is stronger than that for $A_{N-1}$ (\ref{fund weyl chamber A_n}), thus making our life easier.\\
A direct check shows that the thesis holds true for $N=4$.\\
So let us suppose we are in the case $D_N$, with $N \geq 5$. The root $\nu$
is a vector of $\Z^N$ having at most 4 non-zero elements
$a_i \in \{\pm 1, \pm 2\}$ (thus since $N\ge 5$ there is at least one coordinate equal to 0).\\
Look at the last coordinate $a_N$: if $a_N \neq 0$, then there is no way for condition (\ref{fund weyl chamber D_n}) to be satisfied, and we are done. 
If instead $a_N=0$, we look at the other coordinates: if all the other elements are non-zero, then it means that we are necessarily
in the case $D_5$, and such elements are alternating 1's and -1's, hence we are done again. If there is another $a_i=0$ with $i \neq N-1$,
we are also done.\\
Finally, if we are in the case $a_{N-1}=a_N=0$, we repeat the argument above.\\
We can go on until we either get to $a_1=a_2=\ldots=a_N=0$, or we
encounter an element that cannot satisfy  (\ref{fund weyl chamber
D_n}), and this concludes the proof for the $D$-case.

The proof for the three exceptional cases $E_6$, $E_7$ and $E_8$ is
nothing but a brute force check that one can do using any computer
algebra system.
\end{proof}

\begin{thm}\label{teorema tangenti semplici}
Let $\T_{G/P}$ the tangent bundle on a flag manifold $G/P$, where
$G$ is a complex simple Lie group of type $ADE$, and $P$ one of
its parabolic subgroups. Then $\T_{G/P}$ is simple.
\end{thm}
\begin{proof}
Theorem \ref{i tangenti sono 1r} together with Remark \ref{rem sul pullback} imply that the isotypical component $\HH^0(\End(\T_{G/P}))^\C=\HH^0(\End(\TT_{G/P}))^\C=\C$. 
But from Lemma \ref{struttura H0 graduato} we get that $\HH^0(\End(\T_{G/P}))=\HH^0(\End(\T_{G/P}))^\C$, and we are done.
\end{proof}

\section{Stability}

\subsection{Simplicity and stability}
Let us start this section with some basic definitions.
\begin{defn}
Let $H$ be an ample line bundle on a projective variety $X$ of dimension $d$. For any coherent sheaf $E$ on $X$ define the
\emph{slope} $\mu_H(E)$ as:
$$\mu_H(E)=\frac{c_1(E)\cdot H^{d-1}}{\rk E}.$$
$E$ is called \emph{H-stable} (respectively \emph{H-semistable}) if for every coherent subsheaf $F \subset E$ such that $E/F$
is torsion free and $0 < \rk F < \rk E$ we have:
$$\mu_H(F)<\mu_H(E)\:\:\:(\hbox{respectively}\:\:\leq).$$
This notion of stability is known as Mumford-Takemoto stability.
\end{defn}
\begin{defn}
In the same setting as above, $E$ is called \emph{H-polystable} if it decomposes as a direct sum of $H$-stable
vector bundles with the same slope.
\end{defn}

It is a well-known fact that for vector bundles stability implies polystability and the latter implies semistability, 
see for example \cite{koba}. Also stability implies simplicity, and the viceversa is not true in general (see \cite{OSS}, or
\cite{Faini} for a homogeneous counterexample).

We now want to look at our homogeneous vector bundles from the point
of view of differential geometry. A homogeneous variety $G/P$ is in particular a homogeneous K\"{a}hler manifold. 
For an exhaustive introduction on K\"{a}hler-Einstein manifolds we refer the reader to \cite{besse}. 
Here we content ourselves with quoting the results on K\"{a}hler-Einstein and Hermite-Einstein structures
that we need.

The following holds:
\begin{thm}\cite[Theorem 8.95]{besse}\label{esistenza KE}
Every compact, simply connected homogeneous K\"{a}hler manifold admits a unique
(up to homothety) invariant K\"{a}hler-Einstein structure.
\end{thm}

Theorem \ref{esistenza KE} above implies that the tangent bundle $\T_{G/P}$ admits a K\"{a}hler-Einstein structure, and
hence in particular an Hermite-Einstein structure.

If $X$ is a compact K\"{a}hler manifold and $E$ a holomorphic bundle over $X$,
a Hermitian metric on $E$ determines a canonical unitary connection whose curvature
is a $(1,1)$-form $F$ with values in $\End E$.\\
The inner product of $F$ with the K\"{a}hler form is then an endomorphism of $E$.
Metrics which give rise to connections such that the endomorphism is a multiple of the identity are called
Hermite-Einstein metrics.

Indeed, the notion of an Hermite-Einstein connection originated in physics.
Hitchin and Kobayashi made a very precise conjecture connecting this notion to that of Mumford-Takemoto stability,
which is known as the Hitchin-Kobayashi correspondence. Uhlenbeck and Yau showed in \cite{UY} that the
conjecture holds true for compact  K\"{a}hler manifolds.
\begin{thm}\cite{UY}
A holomorphic vector bundle over a compact K\"{a}hler manifold admits an
Hermite-Einstein structure if and only if it is polystable.
\end{thm}

As an immediate consequence we get that:
\begin{cor}
Let $G/P$ a rational homogeneous variety of type $ADE$.
Then the tangent bundle $\T_{G/P}$ is polystable with respect to the anticanonical polarization $-K_{G/P}$ induced by the
Hermite-Einstein structure.
\end{cor}
Recall now that a simple bundle is in particular indecomposable. Thus the polystability of the
tangent bundles $\T_{G/P}$ combined with their simplicity implies that the direct sum of stable bundles in
which they decompose is reduced in reality to only one summand, or in other words that:
\begin{thm}\label{i tangenti sono stabili}
Let $G/P$ a rational homogeneous variety of type $ADE$.
Then the tangent bundle $\T_{G/P}$ is stable with respect to the anticanonical polarization $-K_{G/P}$ induced by the
Hermite-Einstein structure.
\end{thm}

\subsection{Some bounds on stability and polarizations in the $A_n$ case}

A natural question arising from Theorem \ref{i tangenti sono stabili} is whether or not there are
other polarizations having the same property of the anticanonical one, and in case we get a
positive answer, can we describe them? This section contains an answer to these questions in some specific cases: in particular here we assume $G=\SL_{n+1}$.

We start with flag manifolds of type $\F(0,n-1,n)$: these are homogeneous varieties of dimension $2n-1$ and of the form $\SL_{n+1}/P$,
where $P=P(\Sigma)$ is the parabolic obtained removing only the first and the last simple root of the Lie algebra, i.e. $\Sigma=\{\alpha_1,\alpha_{n}\}$.
\begin{prop}\label{caratterizzazione polarizzazione}
Let $\F=\F(0,n-1,n)$ be the flag manifold of the form $\SL_{n+1}/P(\alpha_1,\alpha_n)$, and set:
$$m(n)=\frac{-n+n\sqrt{4n^2+4n-3}}{2(n^2+n-1)}.$$
Then the tangent bundle $\T_\F$ is stable with respect to the polarization $\OO_\F(a,b)$ if and only if it is semistable if and only if:
$$m(n)a \leq b \leq m(n)^{-1}a.$$
\end{prop}
\begin{proof}
Start by noticing that:
$$\F(0,n-1,n)=\PP(Q_{\PP^n}),$$
meaning that we can look at our varieties as the projectivization $\F=\PP(Q_{\PP^n})$
of the quotient bundle $Q_{\PP^n} \simeq \T_{\PP^n}(-1)$ on $\PP^n$. Hence we get projections:
$$\xymatrix{&\F \ar[dr]^\beta\ar[dl]_\alpha&\\
\PP^n&&{\PP^n}^\vee}$$ Hence $\Pic (\F(0,n-1,n))=\Z^2$ is
spanned by
$F=\beta^*\OO_{{\PP^n}^\vee}(1)=\OO_\F(1,0)$ and
$G=\alpha^*\OO_{\PP^n}(1)=\OO_\F(0,1)$.\\
Recalling that the elements of $\F$ are couples $(p,H)$=(point,
hyperplane) such that $p \in H \subset \PP^n$ we also get the
identification:
$$F=\{(p,H)\:|\:p\in H,p=p_0\}\:\:\:\:\hbox{and}\:\:\:\:G=\{(p,H)\:|\:p\in H,H=H_0\}.$$
Moreover, we have two short exact sequences:
\begin{align}
&0 \rightarrow \OO_\F \rightarrow \pi^*Q^* \otimes \OO(1)_{\hbox{rel}} \rightarrow \T_{\hbox{rel}} \rightarrow 0 \label{sequ1}\\
&0 \rightarrow  \T_{\hbox{rel}} \rightarrow \T_\F \rightarrow
\pi^*\T_{\PP^n} \rightarrow 0. \label{sequ2}
\end{align}

All in all, the (quiver associated to the) tangent bundle $\T_\F$ to
these varieties has the simple look:
\begin{equation}\label{quiver tangente proiettivizz}
\xymatrix{\bullet^{n-1}&\\
            \bullet^{1}\ar[r]\ar[u]&\bullet^{n-1}}
\end{equation}
So $\gr \T_\F$ has 3 irreducible summands, all with multiplicity 1
and whose rank is described in the picture (\ref{quiver tangente proiettivizz}) above.\\
Now that we have understood the tangent bundle, we can look at its
subbundles. Of course there are more than two subbundles. Yet it is enough
to check the stability condition only on the two \underline{homogeneous}
subbundles, in virtue of a criterion given by Rohmfeld in his paper \cite{Rohm}, and later
refined by Faini in \cite{Faini}:
\begin{thm}[Rohmfeld-Faini]\label{rohmfeld-faini}
Let $H$ be an ample line bundle. If a homogeneous bundle
$E=E_{\rho}$ is not $H$-semistable then there exists a homogeneous
subbundle $F$ induced by a subrepresentation of $\rho$ such that
$\mu_H(F) > \mu_H(E)$.
\end{thm}

Since our tangent bundles have the particular configuration showed in (\ref{quiver tangente proiettivizz}), all we need to
do is just look at the polarizations $H=aF+bG$ such
that:
$$\left\{
  \begin{array}{l}
\mu_H(E')< \mu_H(\T_\F)\\
\mu_H(E'')<\mu_H(\T_\F)
  \end{array}
\right.$$
where the subbundles $E'$ and $E''$ are the irreducible rank $n-1$
bundles in the following picture:
$$E'=\xymatrix{\bullet &\\
            \circ &\circ}\:\:\:\:\:\:\:\:\:\:\:\:\:\:\:\:\:\:\:\:
E''=\xymatrix{\circ &\\
            \circ &\bullet}
$$

Knowing all the weights of the representation associated to our
bundles, we easily compute their first Chern class, so that the two inequalities above read:
\begin{equation}\label{condizioni sulla slope}
\left\{
  \begin{array}{l}
\frac{(-\lambda_1+n\lambda_2)\cdot(a\lambda_1+b\lambda_2)^{2n-2}}{n-1}
<\frac{(n\lambda_1+n\lambda_2)\cdot(a\lambda_1+b\lambda_2)^{2n-2}}{2n-1}\\
\frac{(n\lambda_1-\lambda_2)\cdot(a\lambda_1+b\lambda_2)^{2n-2}}{n-1}
<\frac{(n\lambda_1+n\lambda_2)\cdot(a\lambda_1+b\lambda_2)^{2n-2}}{2n-1}
  \end{array}
\right.
\end{equation}

Intersection theory is also easy to understand in this particular case; out of all products
$$F^iG^j,\:\:\:\:\:\:\:\:\:\:\:\:\:\:\:i+j=2n-1,$$
the only non-vanishing ones will be (recall that we are pulling back from a $\PP^n$!):
$$\left\{
  \begin{array}{l}
F^{n-1}G^n=F^nG^{n-1}=1\\
F^iG^j=0\:\:\:\:\hbox{for}\:i,j \neq n,n-1,\:i+j=2n-1
  \end{array}
\right.$$

We stress the fact that there is a complete symmetry $\PP^n \leftrightarrow {\PP^n}^\vee$, and thus $F \leftrightarrow G$.\\
Simplifying \ref{condizioni sulla slope} becomes the condition $\T_\F$ is stable if and only if:
$$\left\{
\begin{array}{l}
(n^2+n-1)b^2+nab-n^2a^2>0\\
-n^2b^2+nab+(n^2+n-1)a^2>0
\end{array}
\right.
$$
And from these two inequalities one easily gets that  $\T_\F$ is stable with respect to $H=aF+bG$ if and only $m(n)a \leq b \leq m(n)^{-1}a$, where $m(n)$ is defined as in the statement of the Proposition.

The only thing left to check is the equivalence ``stable $\Leftrightarrow$ semistable'', but this simply follows from the fact that the conditions for
semistability are just the conditions (\ref{condizioni sulla slope}) where we substitute the sign $<$ with a $\leq$.\\
But in reality equality never holds, for $m(n)$ is an irrational coefficient (because $4n^2+4n-3=(2n+1)^2-4$), while we need $(a,b) \in \Z^2$. 
\end{proof}

\begin{rem}
An interesting observation is that as $n$ grows bigger $m(n)$ approaches 1. Hence the cone of polarizations with respect to which $\T_\F$ is stable
collapses to the line $a=b$, which is the once corresponding to the anticanonical. The collapsing process is illustrated in Figure \ref{figura coni} below,
where we have drawn the cone for -respectively- $n=2$ and $n=20$ in the space of polarizations $(a,b)$. The dotted line is the anticanonical polarization $\{a=b\}$.
\begin{figure}[ht]
\centerline{
\mbox{\includegraphics[height=1.8in,width=1.8in]{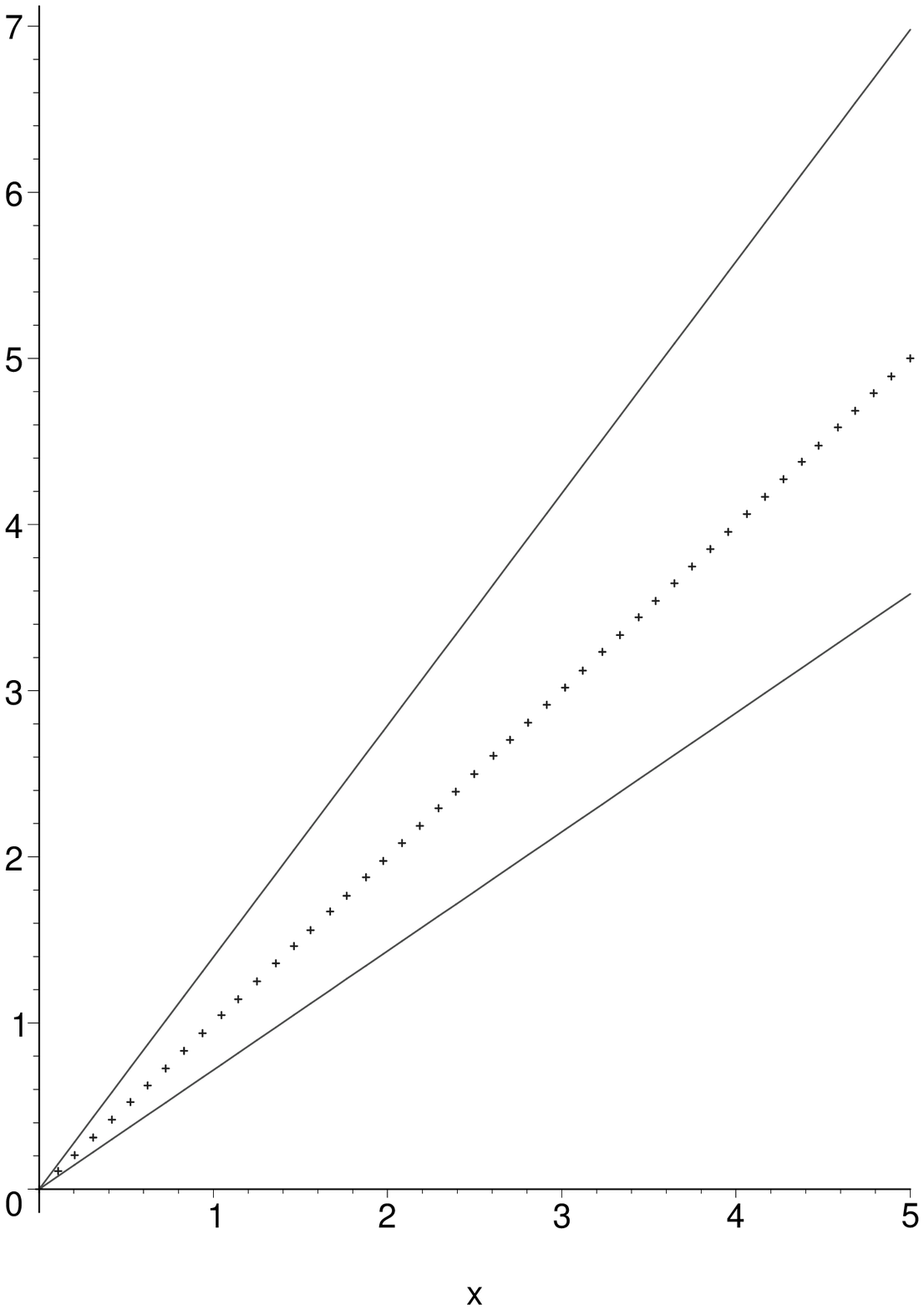}}
\hspace{.35in}
\mbox{\includegraphics[height=1.8in,width=1.8in]{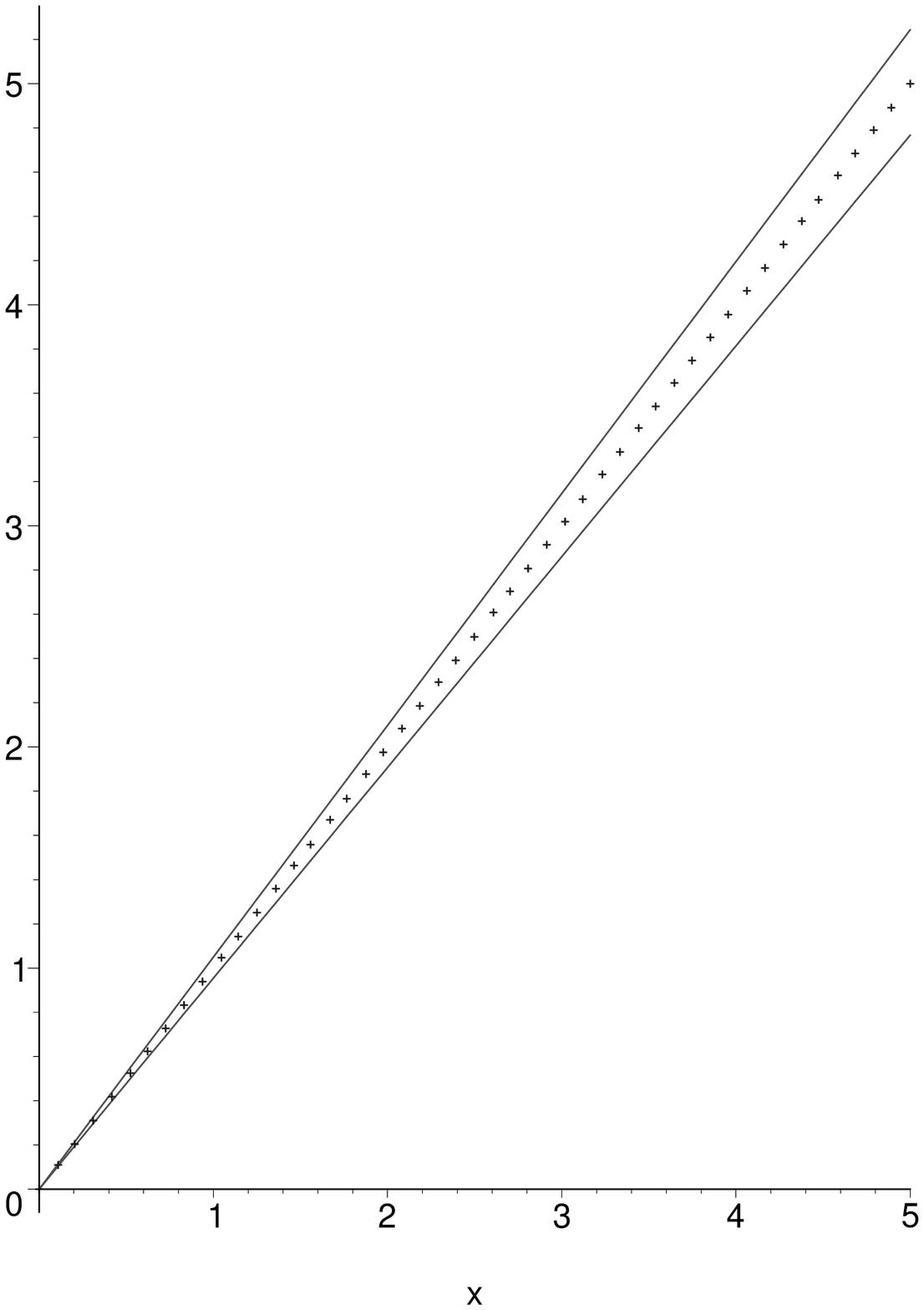}}
}
\caption{The cone for $n=2$ (left) and $n=20$ (right).}
\label{figura coni}
\end{figure}
\end{rem}

We wish to obtain the same type of characterization as Proposition
\ref{caratterizzazione polarizzazione} for other homogeneous varieties.

The next simplest case after the flags with Picard group $\Z^2$ is
the full flag manifold $\F=\SL_4/B$, with $\Pic(\F)=\Z^3$ and dimension 6.\\
The weights of the tangent bundle are the 6 positive roots of $\SL_4$. Before we can show how the the quiver representation $[\T_\F]$ looks like, we need to 
explain here how the relations on the quiver work. In the Borel case, Definition \ref{def relazioni} can be made more explicit.

For each root $\alpha$, let $e_\alpha \in \lieg_\alpha$ be the corresponding Chevalley generator, and define the Chevalley coefficients $N_{\alpha\beta}$ by $[e_\alpha,e_\beta]=N_{\alpha\beta}e_{\alpha+\beta}$, if
$\alpha+\beta \in \Phi^+$, and $N_{\alpha\beta}=0$ otherwise. 

\begin{prop}\cite[Proposition 1.21]{ACGP}\label{def relazioni bis}
The relations $\mathcal{R}$ on the quiver $\Q_{G/P}$ are the ideal
generated by all the equations:
$$\mathcal{R}_{(\alpha,\beta)}=f_{\lambda\mu}f_{\mu\nu}-f_{\lambda\mu'}f_{\mu'\nu}-N_{\alpha\beta}g_{\lambda\nu}$$
for $\alpha < \beta \in \Phi^+$ and for any couple of weights
$\lambda,\nu \in \Q_0$, where $\alpha+\beta=\lambda-\nu$,
$\mu=\lambda+\alpha$ and $\mu'=\lambda+\beta$.
\end{prop}

Let now $G=\SL_N$. For the sake of simplicity we indicate with $e_{ij}$ the element
with weight $L_i-L_j$.
 
Take a couple of roots $\alpha < \beta \in \Phi^+$, $\alpha=L_i-L_j$
and $\beta=L_h-L_k$ (with $i <j$). Note that in this case the non-zero coefficients $N_{\alpha\beta}$ are all
$\pm 1$. The only possibility for the coefficient $N_{\alpha\beta}$
to be non-zero is if either $j=h$ ($\Rightarrow N_{\alpha\beta}=1$),
or $i=k$ ($\Rightarrow N_{\alpha\beta}=-1$). But then from Proposition \ref{def relazioni bis}:
$$\mathcal{R}_{(L_i-L_h,L_h-L_k)}=e_{ih}e_{hk}-e_{hk}e_{ih}-e_{ik}=0$$
$$\mathcal{R}_{(L_i-L_j,L_h-L_i)}=e_{ij}e_{hi}-e_{hi}e_{ij}+e_{hj}=0.$$
In other words the relations tells us nothing newer than $[e_{ih},e_{hk}]=e_{ik}$.\\
Now suppose $j \neq h$ and $i \neq k$, so that $N_{\alpha\beta}=0$.
The relations are thus:
$$\mathcal{R}_{(L_i-L_j,L_h-L_k)}=e_{ij}e_{hk}-e_{hk}e_{ij}=0,$$
for any $i,j,h,k+1,\ldots,N$.

All in all the relations that we need to put on the quiver $\Q_\F$
for a full flag manifold $\F=\SL_N/B$ are nothing but the Serre relations:
\begin{align}
&\mathcal{R}_{(L_i-L_h,L_h-L_k)}=[e_{ih},e_{hk}]=e_{ik},\label{relazioni flag1}\\
&\mathcal{R}_{(L_i-L_j,L_h-L_i)}=[e_{ij},e_{hi}]=-e_{hj},\label{relazioni flag2}\\
&\mathcal{R}_{(L_i-L_j,L_h-L_k)}=e_{ij}e_{hk}-e_{hk}e_{ij}=0,\:\:\hbox{for
} i \neq k,\:j \neq h\label{relazioni flag3}.
\end{align}

For $N=4$, we will have in particular that $[e_{12},e_{34}]=0$, so the corresponding arrows commute.

For $SL_4/B$ the quiver representation $[\T_\F]$ looks like in (\ref{quiver tangente SL4}). We have indicated to which element of $\lien$ 
correspond the arrows.
\begin{equation}\label{quiver tangente SL4}
\T_\F=\xymatrix{\bullet&&\\
\bullet \ar@{-->}[u]^{e_{23}}\ar@{~>}[r]^{e_{12}}&\bullet&\\
\bullet\ar[u]^{e_{34}}\ar@{~>}[r]_{e_{12}} &\bullet \ar[u]_{e_{34}}\ar@{-->}[r]_{e_{23}} &\bullet}
\end{equation}

The relations tell us that the square below is commutative:
$$\xymatrix{\bullet \ar@{~>}[r]^{e_{12}} &\bullet\\
\bullet \ar[u]^{e_{34}} \ar@{~>}[r]_{e_{12}}&\bullet \ar[u]_{e_{34}}}$$

Now let's go back to stability computations. Again by Theorem \ref{rohmfeld-faini} all we need to do is identify all
the homogenous subbundles $F$ of the tangent bundle $\T_\F$, and then impose the stability condition $\mu_H(F)<\mu_H(\T_\F)$;
this will give us necessary and sufficient condition for the polarization $H$ to be such that $\T_\F$ is $H$-stable.

All the homogenous subbundles we need to analyze are thus the following six:
$$E_1\:=\:\xymatrix{\bullet&&\\
\circ &\circ&\\
\circ &\circ &\circ}
\:\:\:\:\:\:\:\:\:\:\:\:\:\:\:\:\:\:\:\:\:\:\:\:\:\:\:\:\:\:\:\:\:\:\:
E_2\:=\:\xymatrix{\circ&&\\
\circ &\bullet&\\
\circ &\circ &\circ}$$

$$E_3\:=\:\xymatrix{\circ&&\\
\circ &\circ&\\
\circ &\circ &\bullet}
\:\:\:\:\:\:\:\:\:\:\:\:\:\:\:\:\:\:\:\:\:\:\:\:\:\:\:\:\:\:\:\:\:\:\:
E_4\:=\:\xymatrix{\bullet&&\\
\bullet \ar@{-->}[u]\ar@{~>}[r]&\bullet&\\
\circ &\circ &\circ}$$

$$E_5\:=\:\xymatrix{\circ&&\\
\circ &\bullet&\\
\circ &\bullet \ar[u]\ar@{-->}[r] &\bullet}
\:\:\:\:\:\:\:\:\:\:\:\:\:\:\:\:\:\:\:\:\:\:\:\:\:\:\:\:\:\:\:\:\:\:\:
E_6\:=\:\xymatrix{\bullet&&\\
\bullet \ar@{-->}[u]\ar@{~>}[r]&\bullet&\\
\circ &\bullet \ar[u]\ar@{-->}[r] &\bullet}$$

Now we need to compute the first Chern class of all these bundles.\\
The elements of $\F$ are triples $(p,\ell,\pi)$=(point, line, plane) such that
$p \in \ell \subset \pi \subset \PP^3$. A basis for the Picard group is given by:
$$\left\{
  \begin{array}{l}
F=\{(p,\ell,\pi)\:|\:p \in \pi_0,\: \pi_0\:\hbox{fixed}\:\}\\
G=\{(p,\ell,\pi)\:|\: \ell \cap \ell_0 \neq \emptyset,\: \ell_0\:\hbox{fixed}\:\}\\
H=\{(p,\ell,\pi)\:|\:p_0 \in \pi,\: p_0\:\hbox{fixed}\: \}
  \end{array}
\right.$$

They correspond to the pull-back (via the standard projection) of the tautological bundle $\OO(1)$ from
respectively $\PP^3$ ($\leftrightarrow F$), $\G(1,3)$ ($\leftrightarrow G$) and ${\PP^3}^\vee$ ($\leftrightarrow H$). 
We underline the symmetry between $F$ and $H$.

Since we are looking for all polarizations
$\OO_\F(a,b,c)=aF+bG+cH$ such that for all $i=1,\ldots,6$:
\begin{equation}\label{condizioni da imporre sui 6 sottofibrati}
\frac{c_1(E_i)(aF+bG+cH)^5}{\rk E_i} < \frac{c_1(\T_\F)(aF+bG+cH)^5}{6},
\end{equation}
we are interested in intersections $F^iG^jH^k$, where $i+j+k=6$.\\
Intersection theory brings us to:

$$\left\{
  \begin{array}{l}
FG^4H=HG^4F=F^2G^2H^2=FG^3H^2=F^2G^3H=2\\
F^3G^2H=FG^2H^3=F^3GH^2=F^2GH^3=1\\
F^iG^jH^k=0\:\:\:\:\hbox{for all other}\:i,j,k\:\hbox{s.t.}\:i+j+k=6
  \end{array}
\right.$$

With the help of a computer algebra system and the intersections above we see that
the six inequalities (\ref{condizioni da imporre sui 6 sottofibrati}) define a cone around the line
$\{a=b=c\}$ that corresponds to the anticanonical polarization $-K_\F=\OO_\F(2,2,2)$.

Figure \ref{figura 6 curve} below shows a section of this cone cut by the plane $\{a+b+c=3\}$ orthogonal to the
``anticanonical line''. It is somewhat unexpected that the region that we obtain is not convex. 
In fact from general theory we learn that the area would be convex if we were considering all possible characters in the definition of stability, 
and not just the ones arising from geometric polarizations given by ampla line bundles like in our case.
In the next section we will explain with some more detail the question of characters, stability and moduli spaces.

\begin{figure}[ht]
\begin{center}
\includegraphics[height=2.5in,width=2.5in]{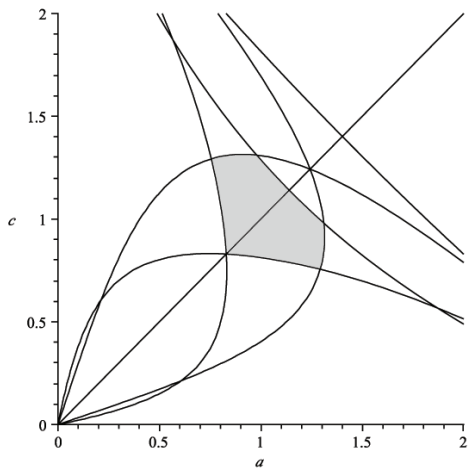}
\caption{Polarizations for $\SL_4/B$}
\label{figura 6 curve}
\end{center}
\end{figure}

\section{Moduli and stability}

There is a notion of semistability of representations of quivers introduced by King in \cite{king},
which is suitable to construct moduli spaces according to the Geometric Invariant Theory (GIT from now on).\\
In their paper \cite{OR} Ottaviani and Rubei showed that King's notion of
semistability for a representation $[E]$ of the quiver $\Q_{G/P}$ is in fact equivalent to the Mumford-Takemoto
semistability of the associated bundle $E$ on $X=G/P$, when the latter is a Hermitian symmetric variety. 
They thus obtain moduli spaces of $G$-homogeneous semistable bundles with fixed $\gr E$.

In this section we recall some of these results and show how they can be extended to our more general setting where $X$
is any -not necessarily Hermitian symmetric- homogeneous variety.\\
Consider the moduli problem of homogeneous vector bundles $E$ on $X$ with the same $\gr E$
and thus with the same dimension vector $\alpha=(\alpha_\lambda) \in \Z^{|\Q_0|}$.

Once we have made the choice of vector spaces $V_\lambda$ with dimension $\alpha_\lambda$, the isomorphism classes
of representations of the quiver $\Q_X$ with the same dimension vector $\alpha$ are in natural 1-1 correspondence with the orbits of the group:
$$\GL(\alpha):=\prod_{\lambda \in \Q_0} \GL(V_\lambda)$$
acting over
$$\mathscr{R}(\Q_X,\alpha):=\oplus_{a \in \Q_1}\Hom(V_{ta},V_{ha})$$
by $(g \cdot \phi)_a=g_{ha}\phi_ag_{ta}^{-1}$, and in particular over the
closed subvariety $V_X(\alpha) \subseteq \mathscr{R}(\Q_X,\alpha)$ defined by the relations in our quiver.\\
The affine quotient $\Spec(\C[V_X(\alpha)]^{\GL(\alpha)})$ is a
single point, represented by $\gr E$ itself, and it thus has no interest for our purposes.

Following \cite{king}, we call a \emph{character} of the category
$\C\Q-mod$ an additive function $\sigma:K_0(\C\Q-mod)\rightarrow \R$ on the Grothendieck group.\\
(For the sake of simplicity we denote by $\C\Q-mod$ the category of left modules on the path algebra $(\C\Q,\mathcal{R})$
of the quiver $\Q$ with relations $\mathcal{R}$: by writing only $\C\Q$ it is understood that we are modding out the
path algebra by the ideal of relations.)

A representation $[E]$ of the quiver is called \emph{$\sigma$-semistable} if $\sigma([E])=0$ and every subrepresentation $[E']\subseteq [E]$
satisfies $\sigma([E'])\leq 0$.\\
Moreover $[E]$ is called \emph{$\sigma$-stable} if the only subrepresentations $[E']$ satisfying $\sigma([E'])=0$ are $[E]$ itself and 0.

When $\sigma$ takes integer values, there is an associated character $\chi_\sigma$ for $\GL(\alpha)$ acting on $\mathscr{R}(\Q_X,\alpha)$. 
More precisely, King shows that the characters of $\GL(\alpha)$ $\chi_\sigma$, $\chi_\sigma:\GL(\alpha)\rightarrow \C^*$ are given by:
$$\chi_\sigma(g)=\prod_{\lambda \in (\Q_X)_0} \det(g_\lambda)^{\sigma_\lambda}$$
for $\sigma \in \Z^{|\Q_0|}$ such that $\sum_\lambda
\sigma_\lambda \alpha_\lambda=0$.

A point in $\mathscr{R}(\Q_X,\alpha)$ corresponding to a representation $[E] \in \C\Q-mod$ is $\chi_\sigma$-semistable
(respectively $\chi_\sigma$-stable) if and only if $[E]$ is $\sigma$-semistable (respectively $\sigma$-stable).

We stress the fact that $\sigma \in \Hom(\C\Q-mod,\Z)$ can be simply seen as an homomorphism that
applied to $E_\lambda$ gives $\sigma_\lambda$.

A function $f \in \C[V_X(\alpha)]$ is called a \emph{relative invariant of
weight $\sigma$} if 
$$f(g \cdot x)=\chi_\sigma(g)f(x),$$ and the space
of such relatively invariant functions is denoted by:
$$\C[V_X(\alpha)]^{\GL(\alpha),\sigma}.$$

So once we have fixed the dimension vector $\alpha$ and a character $\sigma$, we can define the moduli space $M_X(\alpha,\sigma)$ by:
$$M_X(\alpha,\sigma):=\Proj\big(\oplus_{n \geq 0}\C[V_X(\alpha)]^{\GL(\alpha),n\sigma}\big),$$
which is projective over $\Spec(\C[V_X(\alpha)]^{\GL(\alpha)})$,
hence it is a projective variety.

In fact $M_X(\alpha,\sigma)$ has a more geometrical description as the GIT quotient of the open set
$V_X(\alpha)^{ss}$ of $\chi_{\sigma}$-semistable points.

Fix an ample line bundle $H$ (a polarization).\\ 
Every ample line bundle
$H$ defines a character $\sigma_H$ by:
$$\sigma_H(\alpha)_\lambda=\rk(E)(c_1(E_\lambda)H^{n-1})-(c_1(\gr E)H^{n-1})\rk E_\lambda,$$
where $n$ is the dimension of the underlying variety.

Notice that given an $F \in \C\Q-mod$ with dimension vector $\alpha$ and given a fixed character $\sigma$,
we can define the \emph{slope of $F$ with respect to $\sigma$} (or \emph{slope of $\alpha$ w.r.t. $\sigma$}):
$$\mu_\sigma(F)=\mu_\sigma(\alpha)=\sum_{\lambda \in \Q_0}\sigma_\lambda \alpha\lambda.$$
An object is then called $\mu_\sigma$-(semi)stable if and only if it is $\sigma$-(semi)stable.

Recall now from Theorem \ref{theta} that a homogeneous bundle $E$ is determined by
$\theta \in \Hom(\gr E, \gr E \otimes \T_X)$ such that $\theta \wedge \theta = \varphi \theta$.

\begin{thm}\cite[Generalization of Theorem 7.1]{OR}\label{teorema 7.1}
Let $E$ be a homogeneous bundle on a rational homogeneous variety
$X$, and let $\alpha$ be the dimension vector corresponding to $\gr
E$. Fix an ample line bundle $H$ giving the character $\sigma \in \Hom(\C\Q-mod,\Z)$. Then the following facts are equivalent:
\begin{itemize}
  \item[(i)] for every $G$-invariant subbundle $K$, we have $\mu_\sigma(K) \leq
  \mu_\sigma(E)$ \emph{(equivariant semistability)};
  \item[(ii)] for every subbundle $K$ such that $\theta_E(\gr K) \subset \gr K \otimes
  \T_X$, we have $\mu_\sigma(K) \leq  \mu_\sigma(E)$ \emph{(Higgs semistability)};
  \item[(iii)] the representation $[E]$ if $\Q_X$ is
  $\sigma$-semistable, according to \cite{king} \emph{(quiver semistability)};
  \item[(iv)] $E$ is a $\chi_\sigma$-semistable point in $V_X(\alpha)$ for the action of
$\GL(\alpha)$ \cite{king} \emph{(GIT semistability)};
  \item[(v)] for every subsheaf $K$, we have $\mu_H(K) \leq \mu_H(E)$ \emph{(Mumford-Takemoto
  semistability)}.
\end{itemize}
\end{thm}

\begin{proof}
The equivalence (i) $\Leftrightarrow$ (ii) follows from the fact that a subbundle
$K \subset E$ is $G$-invariant if and only if $\theta_E(\gr K) \subset \gr K \otimes \T_X$.
The equivalence (ii) $\Leftrightarrow$ (iii) is just a rephrasing of the second fundamental equivalence of categories.
The equivalence (iii) $\Leftrightarrow$ (iv) is proved in \cite[Proposition 3.1]{king}. In fact this equivalence
holds true even for those characters $\sigma$ that do not have a geometric interpretation as the one induced by a polarization that we chose. 
Finally, the equivalence (i) $\Leftrightarrow$ (v) is proved for example in \cite{migliorini}.
\end{proof}

With the same reasoning one can prove that:

\begin{thm}\cite[Generalization of Theorem 7.2]{OR}\label{teorema 7.2}
Let $E$ be a homogeneous bundle on a rational homogeneous variety
$X$, and let $\alpha$ be the dimension vector corresponding to $\gr
E$. Fix an ample line bundle $H$ giving the character $\sigma \in \Hom(\C\Q-mod,\Z)$. Then the following facts are equivalent:
\begin{itemize}
  \item[(i)] for every proper $G$-invariant subbundle $K$, we have $\mu_\sigma(K) <
  \mu_\sigma(E)$ \emph{(equivariant stability)};
  \item[(ii)] for every proper subbundle $K$ such that $\theta_E(\gr K) \subset \gr K \otimes
  \T_X$, we have $\mu_\sigma(K) <  \mu_\sigma(E)$ \emph{(Higgs stability)};
  \item[(iii)] the representation $[E]$ if $\Q_X$ is
  $\sigma$-stable, according to \cite{king} \emph{(quiver stability)};
  \item[(iv)] $E$ is a $\chi_\sigma$-stable point in $V_X(\alpha)$ for the action of
$\GL(\alpha)$ \cite{king} \emph{(GIT stability)};
  \item[(v)] $E \simeq W \otimes E'$ where $W$ is an irreducible $G$-module, and for every proper
  subsheaf $K \subset E'$, we have $\mu_H(K) < \mu_H(E')$ \emph{(Mumford-Takemoto
  stability)}.
\end{itemize}
\end{thm}

\begin{proof}
The proof is the same of Theorem \ref{teorema 7.1}. The equivalence (i) $\Leftrightarrow$ (v) is proved in \cite{Faini}.
\end{proof}

We remark that case (v) in Theorem \ref{teorema 7.2} is more involved with respect to the naive expectation coming from (v) in 
Theorem \ref{teorema 7.1}.

\bibliographystyle{amsalpha}
\bibliography{grenoble}

\end{document}